\documentclass[journal,9pt]{IEEEtran}

\usepackage{graphics} 
\usepackage{subfigure}
\usepackage{epsfig} 
\usepackage{mathptmx} 
\usepackage{times} 
\usepackage{amsmath} 
\usepackage{amssymb}  
\usepackage{amsthm}
\usepackage{xcolor}
\usepackage{cite}
\usepackage{paralist}

\newtheorem{R}{Remark}
\newtheorem{T}{Theorem}

\newtheorem{Lm}{Lemma}
\newtheorem{Asmp}{Assumption}

\begin{document}
\title{Distributed Subgradient-based Multi-agent Optimization with More
General Step Sizes}
\author{Peng Wang, and Wei Ren\thanks{Peng Wang and Wei Ren are with the Department of Electrical and Computer Engineering,
        University of California, Riverside, Riverside, CA, 92521. Emails:
        {\tt\small pwang033@ucr.edu} and {\tt\small ren@ece.ucr.edu}}}
\date{}
\maketitle
\begin{abstract}
A wider selection of step sizes is explored for the distributed subgradient algorithm for multi-agent optimization problems, for both time-invariant and time-varying communication topologies. The square summable requirement of the step sizes commonly adopted in the literature is removed. The step sizes are only required to be positive, vanishing and non-summable. It is proved that in both unconstrained and constrained optimization problems, the agents' estimates reach consensus and converge to the optimal solution with the more general choice of step sizes. The idea is to show that a weighted average of the agents' estimates approaches the optimal solution, but with different approaches.
In the unconstrained case, the optimal convergence of the weighted average of the agents' estimates is proved by analyzing the distance change from the weighted average to the optimal solution and showing that the weighted average is arbitrarily close to the optimal solution. In the constrained case, this is achieved by analyzing the distance change from the agents' estimates to the optimal solution and utilizing the boundedness of the constraints. Then the optimal convergence of the agents' estimates follows because consensus is reached in both cases.
These results are valid for both a strongly connected time-invariant graph and time-varying balanced graphs that are jointly strongly connected.
\end{abstract}
\section{Introduction}
\label{sec:intr}
With the emergence of large scale networks and complex large systems, distributed optimization arises in many areas such as distributed model predictive control \cite{MotaColor}, distributed signal processing \cite{AmicoSig}, optimal network flow \cite{ZarghamNewton} and network utility maximization \cite{WeiNewton} and has attracted significant attention.
The distributed optimization problems can be roughly classified into two categories.
In the first category, each agent has a local objective function and sometimes a local constraint, both unknown to others, but different agents share the same optimization variable. This means that different agents' estimates of the optimizer should be the same at last \cite{NedicUnconstr,NedicConstr,LobelSubgr,LinSubgr,JokoveticFast,ZhuDualSubgr,ZhuNonconvex}. The problems in this category can be regarded as a distributed potential problem.
In the second category, every agent has a local objective function unknown to others, the constraints of the agents are coupled, and every agent knows only a part of the coupled constraints \cite{ChangPDPert,MotaColor,ZarghamNewton,WeiNewton}. The problems in this category can be regarded as a distributed network flow problem.
In this paper we will focus on the problems in the first category.

Various algorithms have been developed to solve the problems in the first category.
In \cite{NedicUnconstr}, a distributed subgradient algorithm is designed for an unconstrained distributed optimization problem, with the assumption of uniformly bounded subgradients, and a non-degenerate, time-varying and balanced communication topology.
In \cite{NedicConstr}, a distributed optimization problem with identical local constraints or non-identical local constraints in the context of a complete graph is considered through a projected distributed subgradient algorithm. Ref. \cite{LobelSubgr} considers non-identical local constraints for balanced and state-dependent switching graphs.
Then \cite{LinSubgr} proves the convergence of the distributed subgradient algorithm with non-identical local constraints under time-varying balanced and fixed unbalanced graphs.
Some accelerated algorithms are proposed in \cite{JokoveticFast}, in which two distributed Nesterov gradient methods are designed and these algorithms are shown to converge faster than the distributed subgradient algorithm in \cite{NedicUnconstr}.
A zero-gradient-sum algorithm is developed in \cite{LuZGS}, in which each agent starts from its local minimizer and the sum of the gradients is kept at zero.
On the other hand, some dual or primal-dual subgradient algorithms are developed for distributed optimization problems with equality and inequality constraints.
Ref. \cite{ZhuDualSubgr} proposes a distributed primal-dual subgradient algorithm to deal with identical affine equality and convex inequality constraints. A projected subgradient method is designed to find the saddle point of the Lagrangian of the primal problem. Then in \cite{ZhuNonconvex}, a similar idea is adopted to develop a distributed dual subgradient algorithm to solve a non-convex problem approximately, with the consensus requirement relaxed.

In the above papers on subgradient-related distributed solutions to the optimization problem \cite{NedicUnconstr,NedicConstr,LinSubgr,ZhuDualSubgr,ZhuNonconvex}, the step sizes for the subgradient should be positive, vanishing, non-summable but square summable. Intuitively, the positiveness makes the algorithm travel in the descent direction, and the non-summablity makes the subgradient a persistent factor in the optimization process finally leading to the optimal solution. But there seems to be no obvious meaning for the square summability of the step size.

In this paper, we will show that the square summability is not necessary for the distributed subgradient method.
We will prove that a positive, vanishing and non-summable step size can make the agents' estimates converge to the optimal solution in both the unconstrained and constrained distributed optimization problems.
This step size selection is actually the same as that required by the centralized subgradient method \cite{ShorBook}. Our results are valid for both the time-varying balanced and time-invariant unbalanced networks. It is worth mentioning that \cite{LiuContinuous} solves the distributed optimization problem with a continuous-time algorithm, where a feedback term instead of the projection operator is used to drive the agents' estimates to the constraint set. In \cite{LiuContinuous}, the step size is required to be positive and vanishing and to have infinite integral. While the results in \cite{LiuContinuous} are interesting and relax the step size requirement,
our results are different from and complement those in \cite{LiuContinuous} in the following aspects: first, the results in this paper are valid for the problem with non-identical constraints, while those in \cite{LiuContinuous} only deal with that with identical constraints; second, both time-varying balanced graphs and fixed unbalanced graphs are considered in this paper, while only a fixed undirected graph is taken into account in \cite{LiuContinuous}; third, the local objective functions are only required to be convex in this paper, while they are required to be strictly convex and differentiable in \cite{LiuContinuous};
fourth, the algorithms are different (discrete-time algorithm in this paper versus continuous-time algorithm in \cite{LiuContinuous}) and so are the analysis approaches. The discrete-time algorithm is projection based, ensuring that the agents stay in their constraint sets at each time instant while the continuous-time algorithm only ensures that the agents approach their constraint sets eventually.

In this paper, we show that with the more general selection of step sizes, the agents' estimates can still reach a consensus and arrive at an optimal solution using the distributed subgradient method.
For the unconstrained optimization problem, we first show the optimal convergence of a sub-sequence of a weighted average of the estimates of different agents by investigating the distance change from the weighted average to the optimal set. Then we show that as time goes by, the weighted average stays in the neighborhood, vanishing with step size, of arbitrary sublevel sets of the global objective function. Next, with consensus, we prove that the estimates of all agents approach the optimal solution.
For the constrained optimization problem, we first prove the optimal convergence of a sub-sequence of the weighted average of the agents' estimates by studying the distance change from the estimates of different agents to the optimal solution. Then the convergence of the corresponding sub-sequence of the agents' estimates follows from consensus. Next, we show the convergence to the optimal solution of the estimates of different agents with the boundedness of the constraints.
The above results hold for both a strongly connected fixed graph and time-varying balanced graphs that are jointly strongly connected.

\paragraph{Notations}
We use $\mathbb{R}$ for the set of real numbers, $\mathbb{R}^{n}$ for the set of $n\times 1$ real vectors and $\mathbb{R}^{n\times n}$ for the set of $n\times n$ real matrices. The symbol $\mathbb{N}^{+}$ represents the set of positive integers, i.e., $\mathbb{N}^{+}=\{1,2,3,\cdots\}$, and the symbol $\mathbb{N}$ represents the set of natural numbers, i.e. $\mathbb{N}=\{0\}\bigcup \mathbb{N}^{+}$.
A sequence of real numbers or vectors $x(k),\; k=1,2,\cdots,$ is represented by $\{x(k)\}$.
The distance between a point $x$ and some set $X$ is $d(x,X)=\inf\limits_{p\in X}\|x-p\|$, and the distance between two sets $X$ and $Y$ is defined as $d(X,Y)=\inf\limits_{x\in X,\; y\in Y}\|x-y\|$.
 The transpose of a vector $a$ is represented by $a^{T}$. We let $\mathbf{1}_{n}$ be the $n\times 1$ vector of all ones.
We use $P_{X}(x)$ to denote the projection of a point $x$ onto a closed convex set $X$: $P_{X}(x)=\arg\underset{p\in X}\min\|x-p\|$.
The convex hull of a set $X$ is denoted by $\text{conv}(X)$.


\section{Preliminaries}
\label{sec:prel}
In this section, we introduce some preliminary results on graph theory and convex optimization.
\subsection{Graph Theory}
An $n$th order directed graph, denoted by $\mathcal{G}(V,E,A)$, is composed of a vertex set $V=\{1,\cdots,n\}$, an edge set $E\subseteq V\times V$ and an adjacency matrix $A$.
We use the pair $(j,i)$ to denote the edge from vertex $j$ to vertex $i$. We suppose that $(i,i)\in E,\; \forall i\in V$.
The adjacency matrix $A=(a_{ij})_{n\times n}\in \mathbb{R}^{n\times n}$ associated with the graph $\mathcal{G}$ is defined such that $a_{ij}$ is positive if $(j,i)\in E$, and $a_{ij}=0$ otherwise. We assume that $A$ is row stochastic, i.e., $\sum\limits_{j=1}^{n}a_{ij}=1,\;\forall i\in V$.
The graph $\mathcal{G}$ is balanced if $\sum\limits_{j=1}^{n}a_{ij}=\sum\limits_{j=1}^{n}a_{ji},\;\forall i\in V$.
The neighbor set of vertex $i$ is defined as $N_{i}=\{j: (j,i)\in E\}$.
 A directed path from $i$ to $j$ is a sequence of edges $(i,i_{1}),(i_{1},i_{2}),\cdots,(i_{p},j)$, starting from vertex $i$ and sinking at vertex $j$.
 The directed graph $\mathcal{G}$ is strongly connected if for any pair of vertices $i$ and $j$, there is a directed path from $i$ to $j$. Intuitively speaking, every vertex in a strongly connected graph can have some influence on the whole network. The union of a collection of graphs is a graph with the vertex and edge sets being the unions of the vertex and edge sets of the graphs in the collection.
\subsection{Convex Optimization}
A set $C$ is convex if $\forall x,\; y \in C$, $\alpha x+(1-\alpha)y \in C$, $\forall \alpha \in [0,1]$. That is, the line segment is in the set $C$ if the two endpoints are.
The convex hull of a set $D$, denoted by $\text{conv}(D)$ is the smallest convex set that contains $D$, i.e.,
\begin{inparaenum}[\itshape a\upshape)]
 \item $\text{conv}(D)$ is convex,
 \item $D\subset \text{conv}(D)$, and
 \item for arbitrary convex set $C$ that contains $D$, $\text{conv}(D)\subset C$.
\end{inparaenum}
A function $f$ is convex if its domain is convex and for all $x$ and $y$ in its domain,
$f(\alpha x+(1-\alpha)y)\leq \alpha f(x)+(1-\alpha)f(y),\; \forall \alpha \in [0,1]$.

An optimization problem
\begin{align*}
\text{minimize} &\quad f(x)&\quad
\text{subject to} &\quad x\in X
\end{align*}
is a convex optimization problem if the objective function $f(x)$ is convex and the constraint set $X$ is also convex.

A vector $g$ is a subgradient of a function $f$ at the point $x_{0}$ if for all $x$ in the domain of $f$,
\begin{equation}\label{m:subg}
f(x)-f(x_{0})\geq g^{T}(x-x_{0}).
\end{equation}
The set of subgradients of $f$ at $x_{0}$ is called subdifferential, denoted by $\partial f(x_{0})$.
 The concept of subgradients (or subdifferential) is a generalization of that of gradients. When the function $f$ is differentiable at $x_{0}$, the gradient of $f$ at $x_{0}$ is the subgradient.

For a projection operator onto a closed convex set, we have the following non-expansiveness property.
\begin{Lm}\cite{NedicConstr}\label{lm:nexp}
Let $X\subset \mathbb{R}^{m}$ be a closed convex set. For any pair of points $x$ and $y$ in $\mathbb{R}^{m}$, we have $\|P_{X}(x)-P_{X}(y)\|\leq \|x-y\|$.
\end{Lm}

\section{Problem Statement}
\label{sec:prob}
For a multi-agent system with $n$ agents, we regard each agent as a vertex. There is an edge $(j,i)$ if agent $i$ receives information from agent $j$. The corresponding entry $a_{ij}$ in the adjacency matrix $A$ denotes the weight assigned by agent $i$ to the received information from agent $j$.

We will focus on the first kind of distributed optimization problems described in Section \ref{sec:intr}.
Each agent has a private local objective function unknown to the other agents, and shares the same variable with the other agents. Also it has its private local constraint. The goal of the multi-agent system is to cooperatively figure out the minimizer of the weighted sum of all local objective functions in the common part of all local constraints:
\begin{align}\label{m:opt}
\text{minimize} &\quad f(x)=\sum\limits_{i=1}^{n}q_{i}f_{i}(x)&\quad
\text{subject to} &\quad x\in X=\bigcap\limits_{i=1}^{n}X_{i},
\end{align}
where $x$ is the variable of the multi-agent system, $f_{i},\; i\in V,$ are the local objective functions, $X_{i}\subseteq\mathbb{R}^{m}, \; i\in V,$ are the local constraints, and
the positive weights $q_{i},\; i\in V,$ are to be specified later. The problem (\ref{m:opt}) is equivalent to the following problem
\begin{align*}
\text{minimize} &\quad \sum\limits_{i=1}^{n}q_{i}f_{i}(x_{i})
\text{subject to} &\quad x_{i}\in X_{i},\;\forall i\in V,
&\;x_{i}=x_{j},\;\forall i,j\in V,
\end{align*}
where $x_{i}\in \mathbb{R}^{m}$ is the variable of agent $i$.
For an unconstrained problem, we let $X_{i}=\mathbb{R}^{m},\; i\in V$.
Consensus is necessary for this kind of optimization problem, because the variables of different agents should be the same and different agents should figure out a common solution of the problem (\ref{m:opt}).

For the multi-agent network, we have some assumptions on its connectivity and the weights in the adjacency matrices.

\begin{Asmp}\label{connect}
There exists an infinite sequence $k_{0},k_{1},\cdots,k_{p},\cdots$ with $0<k_{p+1}-k_{p}\leq B,\; B\in \mathbb{N}^{+}$, such that the union $\bigcup\limits_{k=k_{p}}^{k_{p+1}-1}\mathcal{G}(k)$ is strongly connected, for all $p\in\mathbb{N}$.
\end{Asmp}
The essence behind Assumption \ref{connect} is that the emerging edges should form a strongly connected graph and these edges should also appear sufficiently often to guarantee consensus and convergence to the optimizer.

\begin{Asmp}\label{stochastic}
The adjacency matrices $A(k),\; k=1,2,\cdots,$ share a common positive left eigenvector associated with eigenvalue $1$. That is,
there exists a constant stochastic vector $q=(q_{1},\cdots, q_{n})^{T}$ with $q_{i}>0,\; i\in V,$ and $\mathbf{1}_{n}^{T}q=1$, such that for all $k$, $q^{T}A(k)=q^{T}$.
\end{Asmp}
\begin{R}
Under Assumption \ref{connect}, when $\mathcal{G}(k)$ is fixed and hence strongly connected, $q$ is the positive left eigenvector of the adjacency matrix $A$ associated with eigenvalue $1$ satisfying $\mathbf{1}_{n}^{T}q=1$.
When the time-varying graph $\mathcal{G}(k)$ is balanced, $q=\frac{1}{n}\mathbf{1}_{n}$.
\end{R}

\begin{Asmp}\label{nondegenerate}
The graph is non-degenerate. That is, there exists $\eta>0$, such that for all $k\in \mathbb{N}$, if $a_{ij}(k)>0$, then $a_{ij}(k)>\eta$, and $a_{ij}(k)=0$ otherwise.
\end{Asmp}
This assumption shows that if agent $i$ receives information from agent $j$, then the edge weight $a_{ij}$ is uniformly bounded away from zero.
This assumption ensures that the influence of an individual agent on the network, if there is any, is persistent and does not vanish as time goes by.

For the optimization problem (\ref{m:opt}), we have the following assumptions:
\begin{Asmp}\label{optset}
The problem (\ref{m:opt}) has a bounded nonempty set of optimal points, denoted by $X^{\star}$.
\end{Asmp}

\begin{Asmp}\label{convex}
Each local objective function $f_{i},\;i\in V,$ is convex and continuous in its local constraint set $X_{i}$.
\end{Asmp}
From \cite{ShorBook} we know that a convex function is continuous in the interior of its domain, but Assumption \ref{convex} only requires $f_{i}$ to be continuous in its local constraint set $X_{i}$.

\begin{Asmp}\label{constr}
Each local constraint set $X_{i},\;i\in V,$ is bounded, closed and convex if $X_{i}\neq \mathbb{R}^{m}$.
\end{Asmp}
As the sum of convex functions is also convex, the global objective function $f$ is convex from Assumption \ref{convex}. With Assumption \ref{constr}, the constraint set $X_{i}$ is convex and so is the intersection $X=\bigcap\limits_{i=1}^{n}X_{i}$.  Then the problem (\ref{m:opt}) is a convex optimization problem.

One of the distributed ways to solve the convex optimization problem (\ref{m:opt}) is to use the distributed subgradient method \cite{NedicUnconstr, NedicConstr, LinSubgr, LobelSubgr}
\begin{align}\label{m:alg}
x_{i}(k+1)=P_{X_{i}}(\sum\limits_{j=1}^{n}a_{ij}(k)x_{j}(k)-\alpha(k)g_{i}(k)),
\end{align}
where $x_{i}(k)$ is agent $i$'s estimate of the minimizer of the global objective function $f$ at the $k$th iteration, $a_{ij}(k)$ is the $(i,j)$th entry of the adjacency matrix $A(k)$ at the $k$th iteration, $\alpha(k)$ is the step size, $g_{i}(k)$ is the subgradient of the local objective function $f_{i}$ at $\sum\limits_{j=1}^{n}a_{ij}(k)x_{j}(k)$, and $P_{X_{i}}$ is the projection operator onto $X_{i}$.

In (\ref{m:alg}), the row stochastic property, i.e., $\sum\limits_{j=1}^{n}a_{ij}=1$, makes each agent reach a consensus and converge to a point minimizing a weighted sum of the local objective functions \cite{LinSubgr}, while the column stochastic property, i.e, $\sum\limits_{j=1}^{n}a_{ji}=1$, makes all agents converge to the optimizer of the sum of the local objective functions \cite{NedicConstr, LinSubgr}.

\begin{Asmp}\label{bound}
The subgradients of $f_{i},\;i\in V,$ are uniformly bounded, i.e., there exists $G>0$ such that for all $g\in\partial f_{i}(x)$, $\|g\|\leq G,\;\forall x\in X_{i},\; \forall i\in V$.
\end{Asmp}
The assumption of uniformly bounded subgradients can be found in many references \cite{NedicUnconstr,NedicConstr,LinSubgr,JokoveticFast,ZhuDualSubgr,ZhuNonconvex}, and plays an important role in the consensus and convergence of the distributed subgradient method. But with Assumption \ref{constr} when the local constraint sets $X_{i},\;\forall i\in V,$ are compact, i.e., closed and bounded, Assumption \ref{bound} is redundant because the boundedness of the subgradients can be deduced from the compactness of the constraint sets.

\begin{Asmp}\label{step}
The step size $\alpha(k)$ is positive, vanishing and non-summable, i.e., $\alpha(k)>0$, $\lim\limits_{k\to\infty}\alpha(k)=0$ and $\sum\limits_{k=0}^{\infty}\alpha(k)=\infty$.
\end{Asmp}
\begin{R}
Assumption \ref{step} allows a wider selection of the step sizes for the distributed subgradient algorithm (\ref{m:alg}), by removing the requirement of $\sum\limits_{k=1}^{\infty}\alpha(k)^{2}<\infty$ commonly adopted in the literature \cite{NedicUnconstr,NedicConstr,LinSubgr,ZhuDualSubgr,ZhuNonconvex}. Also, Assumption \ref{step} is the same as that for the centralized subgradient method \cite{ShorBook}, which might imply that this is among the widest range of step sizes for the distributed subgradient algorithm.
\end{R}

\section{Main Results}
\label{sec:main}
In this section, we prove that all agents' estimates of the minimizer of the convex optimization problem (\ref{m:opt}) generated by the distributed subgradient algorithm (\ref{m:alg}) converge to the optimal solution of (\ref{m:opt}), without requiring $\sum\limits_{k=1}^{\infty}\alpha(k)^{2}<\infty$.
Even without the square summable assumption, the existing results can still ensure that the agents' estimates reach a consensus in both the unconstrained and constrained cases, as summarized in the next lemma.
\begin{Lm}\label{lm:consensus}\cite{NedicConstr,LinSubgr}
For a graph sequence $\mathcal{G}(k),\;k=0, 1,2,\cdots,$ satisfying Assumptions \ref{connect}, \ref{stochastic} and \ref{nondegenerate} and the optimization problem (\ref{m:opt}) satisfying Assumptions \ref{convex}, \ref{bound}, \ref{step} with either $X_{i}=\mathbb{R}^{m},\;i \in V,$ or Assumption \ref{constr}, the agent's estimates $x_{i},\;i\in V,$ in the distributed subgradient algorithm (\ref{m:alg}) reach a consensus, i.e., $\lim\limits_{k\to\infty}\|x_{i}(k)-x_{j}(k)\|=0,\;\forall i,j\in V$.
\end{Lm}

However, it is not clear whether the agents' estimates will converge to the optimal solution.
Next, we will prove the convergence of (\ref{m:alg}) to the optimal solution of (\ref{m:opt}) in both the unconstrained and constrained cases.
\subsection{Unconstrained Case}
In this section, we will prove that the global weighted average of the agents' estimates converges to the global optimal set with the step size in Assumption \ref{step} by analyzing the distance change from the weighted average to the optimal solution. Then as consensus is shown in Lemma \ref{lm:consensus}, all agents reach a common minimizer for the problem (\ref{m:opt}). The rigorous statement is as follows:
\begin{T}\label{Thm:opt}
For a graph sequence $\mathcal{G}(k),\;k=0,1,2,\cdots,$ satisfying Assumptions \ref{connect}, \ref{stochastic} and \ref{nondegenerate} and the optimization problem (\ref{m:opt}) satisfying Assumptions \ref{optset}, \ref{convex}, \ref{bound}, and \ref{step} with $X_{i}=\mathbb{R}^{m},\;i\in V$, the agents' estimates $x_{i},\;i\in V,$ in the distributed subgradient algorithm (\ref{m:alg}) converge to a common point in the optimal set $X^{\star}$ of (\ref{m:opt}).
\end{T}
\begin{it}
Proof:
\end{it}
Let $x^{\star}$ be some point in the optimal set $X^{\star}$. Also let
$y(k)=\sum\limits_{i=1}^{n}q_{i}x_{i}(k)$
be the global weighted average of the estimates of all agents and
\begin{align}\label{m:v}
v_{i}(k)=\sum\limits_{j=1}^{n}a_{ij}(k)x_{j}(k)
\end{align}
be the local weighted average of the estimates of agent $i$'s neighbors.
Then note that $q^{T}A(k)=q^{T}$, i.e. $\sum\limits_{i=1}^{n}q_{i}a_{ij}(k)=q_{j},\;\forall j\in V $ with Assumption \ref{stochastic}, we have
 \begin{align*}
 y(k+1)&=\sum\limits_{i=1}^{n}q_{i}x_{i}(k+1)=\sum\limits_{i=1}^{n}q_{i}(\sum\limits_{j=1}^{n}a_{ij}(k)x_{j}(k)-\alpha(k)g_{i}(k))\\
 &=\sum\limits_{j=1}^{n}(\sum\limits_{i=1}^{n}q_{i}a_{ij}(k))x_{j}(k)-\alpha(k)\sum\limits_{i=1}^{n}q_{i}g_{i}(k)\\
 &=\sum\limits_{j=1}^{n}q_{j}x_{j}(k)-\alpha(k)\sum\limits_{i=1}^{n}q_{i}g_{i}(k)=y(k)-\alpha(k)\sum\limits_{i=1}^{n}q_{i}g_{i}(k).
 \end{align*}
Then the distance between the global weighted average $y(k)$ and the point $x^{\star}$ in the optimal set evolves as follows
$\|y(k+1)-x^{\star}\|^{2}=\|y(k)-\alpha(k)\sum\limits_{j=1}^{n}q_{j}g_{j}(k)-x^{\star}\|^{2}
=\|y(k)-x^{\star}\|^{2}+\alpha(k)^{2}\|\sum\limits_{j=1}^{n}q_{j}g_{j}(k)\|^{2}-2\alpha(k)\sum\limits_{j=1}^{n}q_{j}g_{j}(k)(y(k)-x^{\star}).$
According to Assumption \ref{bound}, the subgradient $g_{j}(k)\leq G$. Note that $g_{j}(k)(y(k)-v_{j}(k))\geq - \|g_{j}(k)\|\|y(k)-v_{j}(k)\|\geq -G\|y(k)-v_{j}(k)\|$ and $f_{j}(v_{j}(k))-f_{j}(y(k))\geq g_{j}^{T}(y(k))(v_{j}(k)-y(k))$ from the definition of subgradients in (\ref{m:subg}), we have
\begin{align*}
&\sum\limits_{j=1}^{n}q_{j}g_{j}(k)(y(k)-x^{\star})\\
=&\sum\limits_{j=1}^{n}q_{j}g_{j}(k)(y(k)-v_{j}(k))+\sum\limits_{j=1}^{n}q_{j}g_{j}(k)(v_{j}(k)-x^{\star})\\
\geq & -\sum\limits_{j=1}^{n}Gq_{j}\|y(k)-v_{j}(k)\|+\sum\limits_{j=1}^{n}q_{j}(f_{j}(v_{j}(k))-f_{j}(x^{\star}))\\
=&-\sum\limits_{j=1}^{n}Gq_{j}\|y(k)-v_{j}(k)\|+\sum\limits_{j=1}^{n}q_{j}(f_{j}(v_{j}(k))-f_{j}(y(k)))\\
&+\sum\limits_{j=1}^{n}q_{j}(f_{j}(y(k))-f_{j}(x^{\star}))\\
\geq & -\sum\limits_{j=1}^{n}Gq_{j}\|y(k)-v_{j}(k)\|+\sum\limits_{j=1}^{n}q_{j}g_{j}(y(k))(v_{j}(k)-y(k))\\
&+\sum\limits_{j=1}^{n}q_{j}(f_{j}(y(k))-f_{j}(x^{\star}))\\
\geq & -2G\sum\limits_{j=1}^{n}q_{j}\|y(k)-v_{j}(k)\|+\sum\limits_{j=1}^{n}q_{j}(f_{j}(y(k))-f_{j}(x^{\star})).
\end{align*}
Combining with the fact that $\|\sum\limits_{j=1}^{n}q_{j}g_{j}(k)\|^{2}\leq \sum\limits_{j=1}^{n}q_{j}\|g_{j}(k)\|^{2}\leq G^{2}$, we have
$\|y(k+1)-x^{\star}\|^{2}
\leq \|y(k)-x^{\star}\|^{2}+\alpha(k)^{2}G^{2}+4\alpha(k)\sum\limits_{j=1}^{n}Gq_{j}\|y(k)-v_{j}(k)\|
-2\alpha(k)\sum\limits_{j=1}^{n}q_{j}(f_{j}(v_{j}(k))-f_{j}(x^{\star})).$

Next we prove that $\liminf\limits_{k\to\infty}\sum\limits_{j=1}^{n}q_{j}(f_{j}(y(k))-f_{j}(x^{\star}))\leq 0$ by contradiction. If not, there exist $\epsilon>0$ and $K_{\epsilon}\in \mathbb{N}^{+}$, such that for all $k>K_{\epsilon}$, $\sum\limits_{j=1}^{n}q_{j}(f_{j}(y(k))-f_{j}(x^{\star}))>\epsilon.$
Then
$\|y(k+1)-x^{\star}\|^{2}
\leq  \|y(k)-x^{\star}\|^{2}+\alpha(k)^{2}G^{2}+4G\sum\limits_{j=1}^{n}q_{j}\|y(k)-v_{j}(k)\|\alpha(k)-2\alpha(k)\epsilon
=\|y(k)-x^{\star}\|^{2}-\alpha(k)\epsilon+G^{2}(\alpha(k)^{2}+\frac{4G\sum\limits_{j=1}^{n}q_{j}\|y(k)-v_{j}(k)\|-\epsilon}{G^{2}}\alpha(k)).$
From Lemma \ref{lm:consensus}, $\lim\limits_{k\to\infty}\|x_{i}(k)-x_{j}(k)\|=0$. We have
\begin{align}\label{m:vtoy}
\begin{split}
&\lim\limits_{k\to\infty}\|v_{i}(k)-y(k)\|\\=&\lim\limits_{k\to\infty}\|\sum\limits_{j=1}^{n}a_{ij}(k)x_{j}(k)-y(k)\|\leq \lim\limits_{k\to\infty}\sum\limits_{j=1}^{n}a_{ij}(k)\|x_{j}(k)-y(k)\|\\
=&\lim\limits_{k\to\infty}\sum\limits_{j=1}^{n}a_{ij}(k)\|x_{j}(k)-\sum\limits_{i=1}^{n}q_{i}x_{i}(k)\|\leq \sum\limits_{j=1}^{n}\sum\limits_{i=1}^{n}\lim\limits_{k\to\infty}\|x_{i}(k)-x_{j}(k)\|\\
=&0.
\end{split}
\end{align}
Then it follows that there exists $K_{c}\in\mathbb{N}^{+}$, such that for all $k>K_{c}$, $\sum\limits_{j=1}^{n}q_{j}\|y(k)-v_{j}(k)\|\leq \frac{\epsilon}{8G}$. Then we have
$\|y(k+1)-x^{\star}\|^{2}\leq \|y(k)-x^{\star}\|^{2}-\alpha(k)\epsilon+G^{2}(\alpha(k)^{2}-\frac{\epsilon}{2G^{2}}\alpha(k)).$
 As $\alpha(k)$ vanishes from Assumption \ref{step}, there exists $K_{\alpha}\in\mathbb{N}^{+}$, such that for all $k>K_{\alpha}$, $\alpha(k)\leq \frac{\epsilon}{2G^{2}}$. Then it follows that $\alpha(k)^{2}-\frac{\epsilon}{2G^{2}}\alpha(k)<0$ and
\begin{equation}\label{m:neg}
\|y(k+1)-x^{\star}\|^{2}\leq \|y(k)-x^{\star}\|^{2}-\alpha(k)\epsilon.
\end{equation}
Denote $K_{0}=\max\{K_{\epsilon},K_{c},K_{\alpha}\}$. We have
$\|y(K_{0}+m)-x^{\star}\|^{2}\leq \|y(K_{0}+1)-x^{\star}\|^{2}-\epsilon\sum\limits_{t=K_{0}+1}^{K_{0}+m-1}\alpha(t)$.
As $\sum\limits_{k=1}^{\infty}\alpha(k)=\infty$, we have $\|y(K_{0}+m)-x^{\star}\|^{2}\leq \|y(K_{0}+1)-x^{\star}\|^{2}-\epsilon\sum\limits_{t=K_{0}+1}^{K_{0}+m-1}\alpha(t)<0$ when $m$ is sufficiently large. This contradicts with the fact that $\|y(K_{0}+m)-x^{\star}\|^{2}\geq 0$.
It follows that
\begin{equation}\label{m:inflim}
\liminf\limits_{k\to\infty}\sum\limits_{j=1}^{n}(f_{j}(y(k))-f_{j}(x^{\star}))\leq 0.
\end{equation}

Next we show the optimal convergence of the agents' estimates. Note that $\sum\limits_{j=1}^{n}q_{j}(f_{j}(y(k))-f_{j}(x^{\star}))\geq 0$ because $x^{\star}$ is in the optimal set $X^{\star}$. Combing with (\ref{m:inflim}), we have that $\liminf\limits_{k\to\infty}\sum\limits_{j=1}^{n}(f_{j}(y(k))-f_{j}(x^{\star}))=0$. Then there exists a sub-sequence $\{y(k_{p})\}$ of $\{y(k)\}$, such that $\lim\limits_{k_{p}\to\infty}y(k_{p})= x^{\star}$ and $\lim\limits_{k_{p}\to\infty}f(y(k_{p}))= f(x^{\star})$, where $f=\sum\limits_{i=1}^{n}q_{i}f_{i}$ as in (\ref{m:opt}).
It follows that for all $\delta>0$, there exists $K_{\delta}\in \mathbb{N}^{+}$, such that for all $k_{p}>K_{\delta}$, $f(y(k_{p}))-f(x^{\star})\leq \delta$.
Define $U_{\delta}=\{y:f(y)-f(x^{\star})=\delta\}$ as the level curve of the global objective function.
Let $d(\delta)=\max\limits_{y\in U_{\delta}}\min\limits_{p\in X^{\star}}\|y-p\|$ be the maximum distance from the level curve $U_{\delta}$ to the optimal set $X^{\star}$. From $\alpha(k)\to 0$ and (\ref{m:vtoy}), there exists $K_{\alpha}^{'}\in\mathbb{N}^{+}$ and $K_{c}^{'}\in \mathbb{N}^{+}$, such that for all $k>K_{\alpha}^{'}$, $\alpha(k)\leq \frac{\delta}{2G^{2}}$ and for all $k>K_{c}^{'}$, $\sum\limits_{j=1}^{n}q_{j}\|y(k)-v_{j}(k)\|\leq \frac{\delta}{8G}$.
If $f(y(k))\leq f(x^{\star})+\delta$, then $\min\limits_{p\in X^{\star}}\|y(k)-p\|\leq d(\delta)$. We have
$\min\limits_{p\in X^{\star}}\|y(k+1)-p\|\leq d(\delta) +\alpha(k)\|\sum\limits_{j=1}^{n}q_{j}g_{j}(k)\|\leq d(\delta)+\alpha(k)G.$
 On the other hand,
if $f(y(k))> f(x^{\star})+\delta$, it follows from (\ref{m:neg}) that when $k>\max\{K_{\alpha}^{'},K_{c}^{'}\}$,
$\|y(k+1)-x^{\star}\|^{2}\leq \|y(k)-x^{\star}\|^{2}-\alpha(k)\delta\leq \|y(k)-x^{\star}\|^{2}.$
Taking into consideration of both cases, we have
when $k>\max\{K_{\alpha}^{'},K_{c}^{'}\}$,
$\min\limits_{p\in X^{\star}}\|y(k+1)-p\|\leq\min\limits_{p\in X^{\star}}\|y(k)-p\|+\max\limits_{k}\{\alpha(k)\}G\leq d(\delta)+\frac{\delta}{2G}.$
As $\delta$ is arbitrary and $d(\delta)\to 0$ when $\delta\to 0$, we get $\min\limits_{p\in X^{\star}}\|y(k)-p\|\to 0$, which means that the global weighted average of all agents' estimates converges to some point in the optimal set $X^{\star}$.
Finally with Lemma \ref{lm:consensus}, we obtain that
$\lim\limits_{k\to\infty}\min\limits_{p\in X^{\star}}\|x_{i}(k)-p\|\leq\lim\limits_{k\to\infty}\min\limits_{p\in X^{\star}}(\|y(k)-p\|+\|x_{i}(k)-y(k)\|)=0,$
which means that the estimates of all agents converge to the optimal set. $\square$

\begin{R}
A similar analysis can be applied to the push-sum subgradient algorithm in \cite{NedicPushsum} to prove that the more general step sizes that are positive, vanishing and non-summable can also guarantee the optimal convergence of the push-sum subgradient algorithm. As there is no significant difference in the proof, we omit it in this paper.
\end{R}

\subsection{Constrained Case}
In this section, we prove that the  distributed subgradient algorithm (\ref{m:alg}) under Assumption \ref{step} without the square summable requirement can drive every agent to the optimal solution of the optimization problem (\ref{m:opt}) with constraints.
\begin{T}\label{Thm:opt_constr}
 For a graph sequence $\mathcal{G}(k),\; k=0,1,2,\cdots,$ satisfying Assumptions \ref{connect}, \ref{stochastic} and \ref{nondegenerate} and the optimization problem (\ref{m:opt}) satisfying Assumptions \ref{optset}, \ref{convex}, \ref{constr}, \ref{bound}, \ref{step}, the agent estimates $x_{i},\;i\in V,$ in the distributed subgradient algorithm (\ref{m:alg}) converge to a common minimizer of (\ref{m:opt}).
\end{T}
\begin{it}
Proof:
\end{it}
Let $x^{\star}$ be some point in the optimal set $X^{\star}$ of the problem (\ref{m:opt}). Let $v_{i}(k)$ be defined in (\ref{m:v}).
Then we have
\begin{align*}
&\sum\limits_{i=1}^{n}q_{i}\|x_{i}(k+1)-x^{\star}\|^{2}\\
=&\sum\limits_{i=1}^{n}q_{i}\|P_{X_{i}}(\sum\limits_{j=1}^{n}a_{ij}(k)x_{j}(k)-\alpha(k)g_{i}(k))-x^{\star}\|^{2}\\
\leq&\sum\limits_{i=1}^{n}q_{i}\|\sum\limits_{j=1}^{n}a_{ij}(k)x_{j}(k)-\alpha(k)g_{i}(k)-x^{\star}\|^{2}\\
=&\sum\limits_{i=1}^{n}q_{i}\|\sum\limits_{j=1}^{n}a_{ij}(k)x_{j}(k)-x^{\star}\|^{2}+\alpha(k)^{2}\sum\limits_{i=1}^{n}q_{i}\|g_{i}(k)\|^{2}\\
&-2\alpha(k)\sum\limits_{j=1}^{n}q_{j}g_{j}^{T}(k)(v_{j}(k)-x^{\star}),
\end{align*}
where the inequality is obtained from Lemma \ref{lm:nexp}.
As $\|\cdot\|^{2}$ is convex and $q^{T}A(k)=q^{T}$, i.e., $\sum\limits_{i=1}^{n}q_{i}a_{ij}(k)=q_{j}$, under Assumption \ref{stochastic}, we have
$\sum\limits_{i=1}^{n}q_{i}\|\sum\limits_{j=1}^{n}a_{ij}(k)x_{j}(k)-x^{\star}\|^{2}\leq \sum\limits_{i=1}^{n}q_{i}\sum\limits_{j=1}^{n}a_{ij}(k)\|x_{j}(k)-x^{\star}\|^{2}
=\sum\limits_{j=1}^{n}(\sum\limits_{i=1}^{n}q_{i}a_{ij}(k))\|x_{j}(k)-x^{\star}\|^{2}
=\sum\limits_{j=1}^{n}q_{j}\|x_{j}(k)-x^{\star}\|^{2}.$
Because $\|g_{i}(k)\|\leq G$ under Assumption \ref{bound}, it follows that
\begin{align}\label{m:x_norm}
\begin{split}
\sum\limits_{i=1}^{n}q_{i}\|x_{i}(k+1)-x^{\star}\|^{2}\leq&\sum\limits_{i=1}^{n}q_{i}\|x_{j}(k)-x^{\star}\|^{2}+\alpha(k)^{2}G^{2}\\
&-2\alpha(k)\sum\limits_{j=1}^{n}q_{j}g_{j}^{T}(k)(v_{j}(k)-x^{\star}).
\end{split}
\end{align}
As $f_{j},\;j\in V,$ are convex,  $ f_{j}(v_{j}(k))-f(x^{\star})\leq g_{j}^{T}(k)(v_{j}(k)-x^{\star})$. As a result, we have
$\sum\limits_{i=1}^{n}q_{i}\|x_{i}(k+1)-x^{\star}\|^{2}
\leq \sum\limits_{i=1}^{n}q_{i}\|x_{i}(k)-x^{\star}\|^{2}+\alpha(k)^{2}G^{2}-2\alpha(k)\sum\limits_{j=1}^{n}q_{j}(f_{j}(v_{j}(k))-f_{j}(x^{\star})).$

Next, we prove that $\liminf\limits_{k\to\infty}\sum\limits_{j=1}^{n}q_{j}(f_{j}(v_{j}(k))-f_{j}(x^{\star}))\leq 0$ by contradiction.
If not, there exist $\epsilon>0$ and $K_{\epsilon}\in \mathbb{N}^{+}$, such that $\forall k>K_{\epsilon}$, $\sum\limits_{j=1}^{n}q_{j}(f_{j}(v_{j}(k))-f_{j}(x^{\star}))>\epsilon$.
Then we have
$\sum\limits_{i=1}^{n}q_{i}\|x_{i}(k+1)-x^{\star}\|^{2}\leq \sum\limits_{i=1}^{n}q_{i}\|x_{i}(k)-x^{\star}\|^{2}+\alpha(k)^{2}G^{2}-2\alpha(k)\epsilon
=\sum\limits_{i=1}^{n}q_{i}\|x_{i}(k)-x^{\star}\|^{2}-\alpha(k)\epsilon+(\alpha(k)^{2}G^{2}-\alpha(k)\epsilon).$
As $\lim\limits_{k\to\infty}\alpha(k)=0$, there exists $K_{\alpha}\in\mathbb{N}^{+}$, such that $\forall k>K_{\alpha}$, $0<\alpha(k)<\frac{\epsilon}{G^{2}}$, which implies that $\alpha(k)^{2}G^{2}-\alpha(k)\epsilon<0$. Hence $\forall k>K=\max(K_{\epsilon},K_{\alpha})$, we have
$\sum\limits_{i=1}^{n}q_{i}\|x_{i}(k+1)-x^{\star}\|^{2}\leq
\sum\limits_{i=1}^{n}q_{i}\|x_{i}(k)-x^{\star}\|^{2}-\alpha(k)\epsilon$.
Because $\sum\limits_{k=1}^{\infty}\alpha(k)=\infty$, it follows that when $k$ is sufficiently large
$\sum\limits_{i=1}^{n}q_{i}\|x_{i}(k+1)-x^{\star}\|^{2}\leq \sum\limits_{i=1}^{n}q_{i}\|x_{i}(K+1)-x^{\star}\|^{2}-\sum\limits_{t=K+1}^{k}\alpha(t)\epsilon<0$.
This contradicts with $\|x_{i}(k+1)-x^{\star}\|^{2}\geq 0$.
It can thus be concluded that $\liminf\limits_{k\to\infty}\sum\limits_{j=1}^{n}q_{j}(f_{j}(v_{j}(k))-f_{j}(x^{\star}))\leq 0$.

Next we show the optimal convergence of the agents' estimates.
Define $y(k)=\frac{1}{n}\sum\limits_{i=1}^{n}P_{X}(x_{i}(k))$. Note that $y(k)\in X$ because $X$ is convex.
 We have
\begin{align*}
\|x_{i}(k)-y(k)\|&=\|x_{i}(k)-\frac{1}{n}\sum\limits_{j=1}^{n}P_{X}(x_{j}(k))\|\leq \frac{1}{n}\sum\limits_{j=1}^{n}\|x_{i}(k)-P_{X}(x_{j}(k))\|\\
&\leq \frac{1}{n}(\sum\limits_{j=1}^{n}\|x_{i}(k)-x_{j}(k)\|+\|x_{j}(k)-P_{X}(x_{j}(k))\|).
\end{align*}
 From Lemma \ref{lm:consensus}, we know that $\lim\limits_{k\to\infty}\|x_{i}(k)-x_{j}(k)\|= 0$. So we have $\lim\limits_{k\to\infty}d(x_{j}(k),X_{i})\leq \lim\limits_{k\to\infty}\|x_{i}(k)-x_{j}(k)\|=0,\;\forall i\in V$.
Then it follows that $\lim\limits_{k\to\infty}\|x_{j}(k)-P_{X}(x_{j}(k))\|=\lim\limits_{k\to\infty}d(x_{j}(k),X)=0$.
Hence we have
\begin{align}\label{m:ytox}
\lim\limits_{k\to\infty}\|x_{i}(k)-y(k)\|=0
\end{align}
 and
$\lim\limits_{k\to\infty}\|v_{i}(k)-y(k)\|=\lim\limits_{k\to\infty}\|\sum\limits_{j=1}^{n}a_{ij}(k)x_{j}(k)-y(k)\|\leq \sum\limits_{j=1}^{n}\lim\limits_{k\to\infty}\|x_{j}(k)-y(k)\|
=0.$
As $f_{i},\; i\in V,$ are convex, and continuous in the constraints under Assumption \ref{convex}, we have
$\liminf\limits_{k\to\infty}\sum\limits_{i=1}^{n}q_{i}(f_{i}(y(k))-f_{i}(x^{\star}))=\liminf\limits_{k\to\infty}\sum\limits_{i=1}^{n}q_{i}(f_{i}(v_{i}(k))-f_{i}(x^{\star}))\leq 0$. Also as $y(k)\in X$, it follows that $\sum\limits_{i=1}^{n}q_{i}(f_{i}(y(k))-f_{i}(x^{\star}))\geq 0$. Then we have $\liminf\limits_{k\to\infty}\sum\limits_{i=1}^{n}q_{i}(f_{i}(y(k))-f_{i}(x^{\star}))=0$.
 Therefore, there exists a sub-sequence $\{y(k_{p})\}$ of $\{y(k)\}$, such that
 \begin{align}\label{m:yleq}
 \begin{split}
 \lim\limits_{k_{p}\to\infty}\sum\limits_{i=1}^{n}q_{i}(f_{i}(y(k_{p}))-f_{i}(x^{\star}))&=\liminf\limits_{k\to\infty}\sum\limits_{i=1}^{n}q_{i}(f_{i}(y(k))-f_{i}(x^{\star}))\\
 &= 0.
 \end{split}
\end{align}
As $\{y(k_{p})\}\in X$ is uniformly bounded from Assumption \ref{constr}, $\{y(k_{p})\}$ has a convergent sub-sequence. Without loss of generality, suppose that the convergent sub-sequence is $\{y(k_{p})\}$ itself, with $y_{\infty}$ being its limit point. We also know that $y_{\infty}\in X^{\star}$ from (\ref{m:yleq}).
Without loss of generality, let $x^{\star}=y_{\infty}$. Then we get from (\ref{m:ytox}) that
\begin{align}\label{m:xsubtoy}
\begin{split}
\lim\limits_{k_{p}\to\infty}\|x_{i}(k_{p})-x^{\star}\|&=\lim\limits_{k_{p}\to\infty}\|x_{i}(k_{p})-y_{\infty}\|\\
&\leq \lim\limits_{k_{p}\to\infty}(\|x_{i}(k_{p})-y(k_{p})\|+\|y(k_{p})-y_{\infty}\|)\\
&=0
\end{split}
\end{align}
and $\lim\limits_{k_{p}\to\infty}\sum\limits_{i=1}^{n}q_{i}f_{i}(x_{i}(k_{p}))= f(x^{\star})$.

We then prove the convergence of the estimates $\{x_{i}(k)\},\forall i\in V$ to $x^{\star}$. For the last term in (\ref{m:x_norm}), we have
$\sum\limits_{i=1}^{n}q_{i}g_{i}^{T}(k)(v_{i}(k)-x^{\star})
\leq \sum\limits_{i=1}^{n} q_{i}\|g_{i}(k)\|\|v_{i}(k)-x^{\star}\|\leq G\sum\limits_{i=1}^{n}q_{i}\|\sum\limits_{j=1}^{n}a_{ij}(k)x_{j}(k)-x^{\star}\|
\leq  G\sum\limits_{j=1}^{n} \sum\limits_{i=1}^{n}q_{i}a_{ij}(k)\|x_{j}(k)-x^{\star}\|=G\sum\limits_{i=1}^{n}q_{i}\|x_{i}(k)-x^{\star}\|,$
where the last equality is obtained from the fact that $q^{T}A(k)=q^{T}$ under Assumption \ref{stochastic}.
Then (\ref{m:x_norm}) can be transformed into
$\sum\limits_{i=1}^{n}q_{i}\|x_{i}(k+1)-x^{\star}\|^{2}\leq \sum\limits_{i=1}^{n}q_{i}\|x_{i}(k)-x^{\star}\|^{2}+\alpha(k)^{2}nG^{2}+2\alpha(k)\sum\limits_{i=1}^{n}Gq_{i}\|x_{i}(k))-x^{\star}\|.$
As $\|x_{i}(k))-x^{\star}\|\leq \|x_{i}(k)\|+\|x^{\star}\|,\;i\in V,$ and both the optimal set $X^{\star}$ and the constraint sets $X_{i},\; i\in V,$ are bounded under Assumptions \ref{optset} and \ref{constr}, $\|x_{i}(k))-x^{\star}\|,\;i\in V,$ are also bounded.
With $\lim\limits_{k\to\infty}\alpha(k)=0$, we have
$\limsup\limits_{k\to\infty}\sum\limits_{i=1}^{n}q_{i}\|x_{i}(k+1)-x^{\star}\|^{2}
\leq \liminf\limits_{k\to\infty} (\sum\limits_{i=1}^{n}q_{i}\|x_{i}(k)-x^{\star}\|^{2}+\alpha(k)^{2}nG^{2}+2\alpha(k)\sum\limits_{i=1}^{n}Gq_{i}\|x_{i}(k))-x^{\star}\|
)
=\liminf\limits_{k\to\infty} \sum\limits_{i=1}^{n}q_{i}\|x_{i}(k)-x^{\star}\|^{2}.$
So $\sum\limits_{i=1}^{n}q_{i}\|x_{i}(k)-x^{\star}\|^{2}$ is convergent. From (\ref{m:xsubtoy}) we have
$\lim\limits_{k\to\infty}\sum\limits_{i=1}^{n}q_{i}\|x_{i}(k)-x^{\star}\|^{2}=\lim\limits_{k_{p}\to\infty}\sum\limits_{i=1}^{n}q_{i}\|x_{i}(k_{p})-x^{\star}\|^{2}=0$.
As
$q_{i}>0,\; i\in V,$ under Assumptions \ref{connect} and \ref{stochastic}, it follows that
$\lim\limits_{k\to\infty}\|x_{i}(k)-x^{\star}\|=0$ and $\lim\limits_{k\to\infty}\sum\limits_{i=1}^{n}q_{i}f_{i}(x_{i}(k))=f(x^{\star})$. $\square$

\begin{R}
We can see from Theorem \ref{Thm:opt} and \ref{Thm:opt_constr} that for both the unconstrained and constrained distributed optimization problems (\ref{m:opt}), the positive, vanishing and non-summable step sizes can guarantee the optimal convergence of the distributed subgradient algorithm (\ref{m:alg}). The square summability $\sum\limits_{k=1}^{\infty}\alpha(k)^{2}<\infty$ is not necessary. When $\mathcal{G}(k)$ is balanced, we can get the minimizer of the sum of the local objective functions. When $\mathcal{G}(k)$ is unbalanced and fixed, we can obtain the minimizer of a weighted average of the local objective functions, with the weights being the elements in the positive left eigenvector of the adjacency matrix associated with eigenvalue $1$.
\end{R}
\begin{R}
Sometimes we can transform an unconstrained distributed optimization problem with unbounded subgradients into one with compact constraints, when the graph is balanced. Suppose every local objective function $f_{i}$ is bounded below, i.e., there exists $B\in \mathbb{R}$, such that $f_{i}(x)\geq B$. Without loss of generality, suppose that $f_{i}(x)\geq 0$ (optimizing $f_{i}(x)$ and $f_{i}(x)-B$ is the same). Also assume that the sublevel set of each local objective function is compact. Then the transformation can be achieved in the following steps.
First, initialize all agents' estimates at the same value, i.e., $x_{i}(0)=x_{j}(0),\forall i, j\in V.$
Then each agent runs a consensus algorithm to compute $f(x(0))=\sum\limits_{i=1}^{n}f_{i}(x_{i}(0))$. Next, denote the sublevel set $S_{i}(a)=\{x:f_{i}(x)\leq f(x(0))+a,\;\forall a\geq 0\}$.
Because $f_{i}(x^{\star})\leq \sum\limits_{i=1}^{n}f_{i}(x^{\star})\leq \sum\limits_{i=1}^{n}f_{i}(x_{i}(0))=f(x(0))\leq f(x(0))+a, \forall a>0$,
it follows that $x^{\star}\in S_{i}(a)$.
Thus, the unconstrained optimization problem becomes a constrained one.
\end{R}
\begin{R}
For the constrained case,
if the step sizes are selected in a non-uniform way as $\alpha_{i}(k)=\alpha(k)(1+\delta_{i}(k))$, where $\alpha_{i}(k)$ is the step size of agent $i$ at iteration step $k$, $\alpha(k)$ satisfies $\alpha(k)>0$, $\lim\limits_{k\to\infty}\alpha(k)=0$ and $\sum\limits_{k=1}^{\infty}\alpha(k)=\infty$, and $\lim\limits_{k\to\infty}\delta_{i}(k)=0$, then we can also prove the convergence to the optimal point of the estimates of the agents generated from the distributed subgradient algorithm with a similar analysis process.
First we analyze the distance change from the agents' estimates to the optimal solution as:
\begin{align}\label{m:dist}
\begin{split}
&\sum\limits_{i=1}^{n}q_{i}\|x_{i}(k+1)-x^{\star}\|^{2}\\
\leq & \sum\limits_{i=1}^{n}q_{i}\|x_{i}(k)-x^{\star}\|^{2}+\alpha(k)^{2}\sum\limits_{j=1}^{n}q_{i}\|g_{i}(k)\|^{2}\\
&+\alpha(k)^{2}\sum\limits_{i=1}^{n}q_{i}\|\delta_{i}(k)g_{i}(k)\|^{2}+2\alpha(k)^{2}\sum\limits_{i=1}^{n}q_{i}\delta_{i}(k)g_{i}^{T}(k)g_{i}(k)\\
&-2\alpha(k)\sum\limits_{j=1}^{n}q_{j}(f_{j}(v_{j}(k))-f(x^{\star}))\\
&+2\alpha(k)\sum\limits_{i=1}^{n}q_{i}\delta_{i}(k)\|g_{i}(k)\|\|(v_{i}(k)-x^{\star})\|.
\end{split}
\end{align}
Then we can prove that $\liminf\limits_{k\to\infty}\sum\limits_{j=1}^{n}q_{j}(f_{j}(v_{j}(k))-f_{j}(x^{\star}))\leq 0$ by contradiction. Suppose not. Then there exist $\epsilon>0$ and $K\in \mathbb{N}^{+}$, such that for all $k>K$, $\sum\limits_{j=1}^{n}q_{j}(f_{j}(v_{j}(k))-f_{j}(x^{\star}))>\epsilon$.
Because the local constraint set $X_{i},i\in V,$ is bounded from Assumption \ref{constr}, we have that $v_{i}(k)\in \text{conv}(\bigcup\limits_{i=1}^{n}X_{i})$ and $x^{\star}\in\bigcap\limits_{i=1}^{n}X_{i}$ are bounded. Then it follows that $v_{i}(k)-x^{\star}$ is also bounded, i.e., there exists $H>0$, such that $\|v_{i}(k)-x^{\star}\|\leq H$.
Then we have from (\ref{m:dist}) that
$\sum\limits_{i=1}^{n}q_{i}\|x_{i}(k+1)-x^{\star}\|^{2}
\leq \sum\limits_{i=1}^{n}q_{i}\|x_{i}(k)-x^{\star}\|^{2}+G^{2}(\sum\limits_{i=1}^{n}(1+q_{i}\delta_{i}^{2}+2q_{i}\delta_{i}))\alpha(k)^{2}+2(GH\sum\limits_{i=1}^{n}q_{i}\delta_{i}-\epsilon)\alpha(k)
=\sum\limits_{i=1}^{n}q_{i}\|x_{i}(k)-x^{\star}\|^{2}-\epsilon\alpha(k)+G^{2}(\sum\limits_{i=1}^{n}(1+q_{i}\delta_{i}^{2}+2q_{i}\delta_{i}))\alpha(k)^{2}+(2GH\sum\limits_{i=1}^{n}q_{i}\delta_{i}-\epsilon)\alpha(k),$
As $\delta_{i}\to 0$, there exists $K_{1}\in\mathbb{N}^{+}$, such that for all $k>K_{1}$, we have $q_{i}\delta_{i}^{2}+2q_{i}\delta_{i}<1,\forall i\in V$, and $\sum\limits_{i=1}^{n}q_{i}\delta_{i}\leq \frac{\epsilon}{4GH}$. Then it follows that
$\sum\limits_{i=1}^{n}q_{i}\|x_{i}(k+1)-x^{\star}\|^{2}=\sum\limits_{i=1}^{n}q_{i}\|x_{i}(k)-x^{\star}\|^{2}-\epsilon\alpha(k)+2nG^{2}\alpha(k)^{2}-\frac{\epsilon}{2}\alpha(k).$
As $\alpha(k)\to 0$, there exists $K_{2}\in\mathbb{N}^{+}$, such that for all $k>K_{2}$, we have $\alpha(k)\leq \frac{\epsilon}{4nG^{2}}$. Then it follows that
$\sum\limits_{i=1}^{n}\|x_{i}(k+1)-x^{\star}\|^{2}=\sum\limits_{i=1}^{n}\|x_{i}(k)-x^{\star}\|^{2}-\epsilon\alpha(k)$.
The rest part is very similar to the proof of Theorem \ref{Thm:opt_constr} and is omitted.
\end{R}

\section{Conclusions}
\label{sec:conclu}
We proved the convergence to a common optimal solution of the distributed subgradient method for a distributed convex optimization problem for both the unconstrained and constrained cases. We relaxed the requirement on the step size by removing the square summable requirement, and showed the positive, vanishing and non-summable step sizes were sufficient for the optimal convergence of the distributed subgradient algorithm, when the topology is fixed or time-varying but balanced.

\bibliographystyle{IEEEtran}
\bibliography{IEEEabrv,mybib}
\end{document}